\documentclass[11pt]{article}
\usepackage{theorem}
\usepackage{a4}
\usepackage{amssymb}
\usepackage{amsmath}

{\theorembodyfont{\slshape}
\newtheorem{theorem}{Theorem}[section]
\newtheorem{lemma}[theorem]{Lemma}
\newtheorem{proposition}[theorem]{Proposition}
\newtheorem{corollary}[theorem]{Corollary}

}
{\theorembodyfont{\upshape}

\newtheorem{remark}[theorem]{Remark}

\newtheorem{notit}[theorem]{{}}

}

\newenvironment{proof}[1]{{\em Proof{#1}}.}{\qed}

\newcommand{\qed}{\hspace*{\fill}$\square$\vspace{2ex}}

\newcommand{\beeq}[1]{\begin{eqnarray}\label{#1}}
\newcommand{\eneq}{\end{eqnarray}}

\newcommand{\ke}{{\cal E}}

\newcommand{\kh}{{\cal H}}

\newcommand{\kk}{{\cal K}}
\newcommand{\kl}{{\cal L}}

\newcommand{\ko}{{\cal O}}

\newcommand{\IC}{{\mathbb C}}

\newcommand{\IN}{{\mathbb N}} 
\newcommand{\IP}{{\mathbb P}} 
\newcommand{\IQ}{{\mathbb Q}}

\newcommand{\IZ}{{\mathbb Z}}

\newcommand{\gothS}{{\mathfrak S}}

\newcommand{\dual}{\makebox[0mm]{}^{{\scriptstyle\vee}}}

\newcommand{\isom}{\cong}

\newcommand{\tensor}{\otimes}

\newcommand{\Pic}{{\rm Pic}}

\newcommand{\Hom}{{\rm Hom}}
\newcommand{\Ext}{{\rm Ext}}

\renewcommand{\det}{{\rm det}}

\newcommand{\id}{{\rm id}}

\newcommand{\rk}{{\rm rk}}

\newcommand{\kext}{\ke xt}
\newcommand{\khom}{\kh om}

\newcommand{\lra}{\longrightarrow}
\newcommand{\lla}{\longleftarrow}
\newcommand{\ses}[3]{0\rightarrow#1\rightarrow#2
   \rightarrow#3\rightarrow0}

\newcommand{\verylongarrow}[1]{\raisebox{0mm}[0.9ex][0ex]
{\hbox to #1{\rightarrowfill}}}

\newcommand{\CP}{\IC\IP}

\newcommand{\Ptilde}{\widetilde P}

\begin{document}


\title{On the Cobordism Class of the Hilbert Scheme\\
of a Surface}

\author{Geir Ellingsrud, Lothar G\"ottsche, Manfred Lehn}
 
\date{April 7, 1999}

\maketitle 

\begin{abstract}
\noindent
Let $S$ be a smooth projective surface and $S^{[n]}$ the Hilbert
scheme of zero-dimensional subschemes of $S$ of length $n$. We proof that
the class of $S^{[n]}$ in the complex cobordism ring depends only on the 
class of the surface itself. Moreover, we compute the cohomology and 
holomorphic Euler characteristics of certain tautological sheaves on $S^{[n]}$
and prove results on the general structure of certain integrals over polynomials
in Chern classes of tautological sheaves.
\end{abstract}


Let $S$ be a smooth projective surface over the field of complex numbers. For
a nonnegative integer $n$ let $S^{[n]}$ denote the Hilbert scheme parameterizing
zero-dimensional subschemes of length $n$. By a well-known result of Fogarty 
\cite{Fogarty} 
the scheme $S^{[n]}$ is smooth and projective of dimension $2n$,
and is irreducible if $S$ is irreducible.

Let $\Omega=\Omega^U\tensor\IQ$ be the complex cobordism ring with rational
coefficients. Milnor \cite{Milnor}
showed that $\Omega$ is a polynomial ring
freely generated by the cobordism classes $[\CP_i]$ for $i\in \IN$.
For a smooth and projective complex surface we define 
$$H(S):=\sum_{n=0}^\infty[S^{[n]}]\,z^n,$$
which is an invertible element in the formal power series ring $\Omega[[z]]$.

The main result of this note is the following

\begin{theorem}\label{MainTheorem}--- $H(S)$ depends only on the cobordism 
class
$[S]\in \Omega_2$.
\end{theorem}

>From this we have an immediate corollary:

\begin{corollary}--- If the class of a surface $S$ satisfies the linear
relation $[S]=a_1[S_1]+a_2[S_2]$ for two surfaces $S_1$ and $S_2$ and two
rational numbers $a_1$ and $a_2$, then
$$H(S)=H(S_1)^{a_1}H(S_2)^{a_2}.$$
\end{corollary}

\begin{proof}{} Assume first that $a_1=a_2=1$. The class $[S_1]+[S_2]$ is 
represented by
the disjoint union $S_1\sqcup S_2$. It is clear that the Hilbert scheme of 
a disjoint union satisfies
\beeq{disjointunion}
(S_1\sqcup S_2)^{[n]}=\coprod_{n_1+n_2=n}S_1^{[n_1]}\times S_2^{[n_2]}.
\eneq
Hence we get $H(S)=H(S_1\sqcup S_2)=H(S_1)H(S_2)$. By induction, we get
$H(S)^n=H(S_1)^{n_1}H(S_2)^{n_2}$ whenever there is a relation
 $n[S]=n_1[S_1]+n_2[S_2]$ for positive integers $n_1$, $n_2$ and $n$. The 
corollary follows formally from this.\end{proof}

If $\phi$ is any genus, i.e.\ a ring homomorphism from $\Omega$ to another ring
\cite{Hirzebruch},
this gives $\phi(H(S))=\phi(H(S_1))^{a_1}\phi(H(S_2))^{a_2}$ whenever
$[S]=a_1[S_1]+a_2[S_2]$. Hence if the value of a genus on $H(S)$ is known
for two independent surfaces, the value on $H(S)$ is determined for any surface
$S$. Using this we give a new proof of the following theorem first proved in
\cite{GoettscheSoergel}.

\begin{theorem}\label{ChiYgenus}---
$$\chi_{-y}(H(S))=\exp\left(\sum_{m=1}^{\infty}
        \frac{\chi_{-y^{m}(S)}}{1-(yz)^{m}}\frac{z^{m}}{m}\right).$$
\end{theorem}

\begin{proof}{} Both sides of this equality are multiplicative in $[S]$,
so it suffices to check it for $\CP_2$ and $\CP_1\times\CP_1$.
Now the Hilbert schemes of both $\CP_{2}$ and $\CP_1\times\CP_1$
have a $\IC^{*}$-action with isolated fix-points. For such varieties the
$\chi_{y}$-genus is given by
\begin{displaymath}
        \chi_{-y}=\sum_{x}y^{\dim T(x)^+}=\sum_{p\geq0}^{}b_{2p}y^p,
\end{displaymath}
where the first sum is over all fix-points and $T(x)^{+}$ is the subspace of 
the tangent space at $x$ generated by eigenvectors for the $\IC^{*}$-action
whose eigenvalues have positive weight. In the second sum the $b_{j}$ are the
Betti numbers of the variety. Recall from \cite{GE-SAS} that the odd Betti 
numbers for $\CP_2^{[n]}$ and $(\CP_1\times\CP_1)^{[n]}$ vanish whereas the
even ones are given by
\begin{displaymath}
b_{2p}=\sum_{n_{1}+n_{2}+n_{3}=n}\,\,\sum_{r_3-r_1=p-n}
   p(n_{1},r_{1})p(n_{2})p(n_{3},r_{3}),
\end{displaymath}
and
\begin{displaymath}
b_{2p}=\sum_{n_1+n_2+n_3+n_4=n}\,\,\sum_{r_4-r_1=p-n}
   p(n_1,r_1)p(n_2)p(n_3)p(n_4,r_4),
\end{displaymath}
where $p(a,b)$ denotes the number of partitions of the nonnegative integer
$a$ into $b$ parts. This gives
\begin{displaymath}
\chi_{-y}(\CP_{2}^{[n]})
  =\sum_{p}\sum_{n_1+n_2+n_3=n}\,\,\sum_{r_3-r_1=p-n}
  p(n_1,r_1)p(n_2)p(n_3,r_3)y^{p}.
\end{displaymath}
and
\begin{eqnarray*}
\lefteqn{\chi_{-y}((\CP_1\times\CP_1)^{[n]})}\hspace{4em}\\
&=&\sum_{p}
   \sum_{{n_1+n_2+n_3+n_4=n}}
   \sum_{{r_4-r_1=p-n}}
      p(n_1,r_1)p(n_2)p(n_3)p(n_4,r_4)y^{p}.
\end{eqnarray*}
Now easy calculations using that for any integer $\varepsilon$ we have
$$
\prod_{k\ge 1}^{}\left(1-y^{k+\varepsilon} z^{k}\right)^{-1}
=\sum_{n,r}p(n,r)y^{n+\varepsilon r}z^{n}
$$
show that 
\begin{displaymath}
\chi_{-y}(H(\CP_{2}))=\prod_{k\ge1}
   \left(1-y^{k-1}z^{k}\right)^{-1}
   \left(1-y^{k}z^{k}\right)^{-1}
   \left(1-y^{k+1}z^{k}\right)^{-1}
\end{displaymath}
and
\begin{displaymath}
\chi_{-y}(H(\CP_{1}\times\CP_{1}))=\prod_{k\ge 1}
   \left(1-y^{k-1}z^{k}\right)^{-1}
   \left(1-y^{k}z^{k}\right)^{-2}
   \left(1-y^{k+1}z^{k}\right)^{-1}.
\end{displaymath}
On the other hand we have $\chi_{-y}(\CP_{2})=1+y+y^2$
and $\chi_{-y}(\CP_{1}\times\CP_{1})=1+2y+y^2$, and a standard calculation shows
that these $\chi_y$--genera are related by the exponential expression in the 
theorem.
\end{proof}

For  integers $N$  and $k$ with $0\le k\le N$ there is a genus
$\phi_{N,k }$ whose characteristic power series is
\begin{displaymath}
        x\frac{e^{-\frac{k}{N}x}}{1-e^{-x}}.
\end{displaymath}
If $X$ is a variety whose canonical line bundle has an $N$-th root $L$,
the genus $\phi_{N,k}$ has a geometric interpretation as
$\chi(X,L^{\otimes k})$. For all varieties $\phi$ is
the value of the level-$N$ elliptic genus of Hirzebruch
\cite{HirzebruchBergerJung} in one of the cusps. We have the following

\begin{theorem}\label{phiNk}--- 
        $$\phi_{N,k}(H(S))=\frac{1}{(1-t)^{\phi_{N,k}(S)}}.$$
\end{theorem}

\begin{proof}{} For a line bundle $L$ on $S$ let
$L_n:=f^*g_{*}(\tensor_{i=1}^{n}pr_{i}^{*}L)^{\gothS_{n}}$,
where $f:S^{[n]}\to S^{(n)}$ and $g:S^{n}\to S^{(n)}$ are the two natural
maps and where $pr_{i}$ is the projection of $S^n$ onto the $i$-th factor.
We will  show later that $\chi(L_n)=\binom{\chi(L)+n-1}{n}$, cf.\ 
Lemma \ref{chivonLn}.
Assume that $S$ is a surface such that $\omega_{S}=L^{\otimes N}$
for some line bundle $L$. Then we claim that the same is true for
the Hilbert scheme $S^{[n]}$. Indeed, we have $\omega_{S^{[n]}}=(\omega_S)_n$,
and $L_{n} $ is an $N$-th root of $\omega_{S^{[n]}}$.
This shows the theorem for surfaces having an $N$-th root of the canonical
line bundle. The formula in the theorem being multiplicative, it will be
sufficient to find two independent surfaces having this property. Indeed,
a K3 surface and a product of two curves of the same genus $g$ such that
$N\vert (2g-2)$ will do.
\end{proof}

The strategy for proving Theorem \ref{MainTheorem} is this: First recall 
that the cobordism class of a stably complex manifold is completely determined
by the collection of its Chern numbers. Thus the theorem is equivalent to the
following proposition.

\begin{proposition}\label{propB}--- For any integer $n$ and any partition
 $\lambda$ of $2n$ 
there is a universal polynomial $P_\lambda\in \IQ[z_1,z_2]$ such that 
 the following relation holds for every smooth projective surface $S$:
$$c_\lambda(S^{[n]})=P_\lambda(c_1^2(S),c_2(S)).$$
\end{proposition}

This proposition will be proved by induction on $n$.

The second named author would like to thank Daniel Huybrechts, who explained 
to him the natural quadratic form on the second cohomology of complex 
symplectic manifolds and how it can be used to obtain the holomorphic
Euler characteristic of line bundles on the Hilbert scheme of points on
K3-surfaces. We thank M.S.\ Narasimhan for very useful discussions.


\section{The geometric set-up}


Let $\Sigma_n\subset S^{[n]}\times S$ be the universal family of subschemes
parameterized by $S^{[n]}$, and let $I_n\subset \ko_{S^{[n]}\times S}$ and 
$\ko_n$ denote its ideal sheaf and structure sheaf, respectively.
\beeq{dispA}
\begin{array}{ccccc}
\Sigma_n&\subset&S^{[n]}\times S&\stackrel{q}{\lra}&S\\
&&\scriptstyle{p}\Big\downarrow\phantom{\scriptstyle{p}}\\
&&S^{[n]}.
\end{array}
\eneq
The induction step involves the incidence variety $S^{[n,n+1]}$ of all pairs
$(\xi,\xi')\in S^{[n]}\times S^{[n+1]}$ satisfying $\xi\subset \xi'$. If $\xi'$
is obtained by extending $\xi$ at the closed point $x\in S$, there is an 
exact sequence
\beeq{version1}
0\lra{I_{\xi'}}\lra{I_\xi}\stackrel{\lambda}{\lra}k(x)\lra0.
\eneq
Observe that $\lambda$ defines a point in the fibre of the morphism
$$\sigma=(\phi,\rho):\IP(I_n)\lra S^{[n]}\times S$$
over $(\xi,x)$. (Here for any coherent sheaf $F$ we denote
 $\IP(F):=Proj(Sym^*(F))$, and $\ko_F(1)$ is the tautological quotient
line bundle.) Conversely, any point in the fibre of $\sigma$ defines an 
extension $\xi'$.
In fact, let $\kl:=\ko_{I_n}(1)$, let $\Gamma\subset \IP(I_n)\times S$ 
be the graph of $\rho$, and let $j:=(\id,\rho)$. Then the kernel of the 
composite epimorphism
$$\beta:(\phi\times\id_S)^*I_n\lra (\phi\times\id_S)^*I_n|_{\Gamma}
\isom j_*\sigma^*I_n\lra j_*\kl$$
is an $\IP(I_n)$-flat family of ideal sheaves and induces a classifying 
morphism
$$\psi:\IP(I_n)\lra S^{[n+1]},$$
such that $Ker(\beta)=\psi_S^*(I_{n+1})$. 
This leads to a scheme theoretic isomorphism $S^{[n,n+1]}\isom\IP(I_n)$ and
the basic diagram
\beeq{dispB}
\begin{array}{ccccc}
S&\stackrel{\rho}{\lla}&\IP(I_n)&\stackrel{\psi}{\lra}&S^{[n+1]}\\[1ex]
&&\phantom{\scriptstyle{\phi}}\Big\downarrow\scriptstyle{\phi}\\[2ex]
&&S^{[n]}
\end{array}
\eneq
together with two short exact sequences of universal families
\beeq{basicseq}\ses{\psi_S^*I_{n+1}}{\phi_S^*I_n}{j_*\kl}
\eneq
\beeq{basicseq2}\ses{j_*\kl}{\psi_S^*\ko_{n+1}}{\phi_S^*\ko_n}.
\eneq
Note that $j_*\kl=p^*\kl \otimes\ko_\Sigma=p^*\kl\otimes \rho_S^*\ko_\Delta$.
Here and throughout the paper we will use the short form $f_S:=f\times \id_S$
in order to simplify the notations. 

Since $S$ is smooth of dimension $2$ and $I_n$ is $S^{[n]}$-flat and fibrewise
torsion free, there is a resolution of length 1
$$0\lra{A}\lra{B}\lra{I_n}\lra 0$$
by locally free $\ko_{S^{[n]}\times S}$-sheaves $A$ and $B$ of
rank $a$ and $a+1$, respectively. Then $\IP(I_n)$ embeds into $\IP(B)$ as
the zero locus of the 
homomorphism
$$\pi^*A\lra\pi^*B\lra\ko_B(1)$$
(where $\pi:\IP(B)\to S^{[n]}\times S$ is the projection) such that
 $\kl=\ko_B(1)|_{\IP(I_n)}$. The incidence variety $S^{[n,n+1]}$ is irreducible
\cite{EllStrIrr, CheahSmooth}. (In fact, it is also smooth, though we will not
use this. This result has been independently proven by Cheah 
\cite{CheahSmooth}, Ellingsrud (unpublished) and Tikhomirov \cite{TikhoSmooth}.)
As $S^{[n,n+1]}$ has the correct dimension, 
it is a locally complete intersection and its fundamental class is given by
$$[\IP(I_n)]=c_a(\kl\tensor\pi^*A\dual)\in A_{2n+2}(\IP(B)).$$

\begin{lemma}--- Let $\ell=c_1(\kl)$. Then for any class
$u\in A_*(S^{[n]}\times S)$ one has
\beeq{pushdownofell}
\sigma_*(\ell^i)=(-1)^ic_i(\ko_n)=(-1)^ic_i(-I_n).
\eneq
\end{lemma}

\begin{proof}{} Let $\varepsilon:=c_1(\ko_B(1))$. Then
\begin{eqnarray*}
\sigma_*(\ell^i)&=&\pi_*(\varepsilon^i\cdot[\IP(I_n)])=\pi_*(\varepsilon^i 
c_a(\pi^*A\dual\tensor\ko_B(1))\\
&=&\pi_*(\sum_{k=0}^a\varepsilon^{i+a-k}\pi^*c_k(A\dual))\\
&=&\sum_{k=0}^as_{i-k}(B\dual)c_k(A\dual)\\
&=&c_i(A\dual-B\dual)=(-1)^ic_i(\ko_n).
\end{eqnarray*}
\end{proof}

\begin{proposition}\label{hilbertchowvanishing}--- Let $f_n:S^{[n]}\lra S^{(n)}$ denote the Hilbert-Chow morphism from the Hilbert scheme to the symmetric product. Then
$$(f_n)_*\ko_{S^{[n]}}=\ko_{S^{(n)}}\quad\mbox{and}\quad R^i(f_n)_*\ko_{S^{[n]}}=0\quad\mbox{for all }i>0.$$
\end{proposition}

\begin{proof}{}
 As $f_n$ is a birational proper morphism of normal varieties, 
it follows that $(f_n)_*\ko_{S^{[n]}}=\ko_{S^{(n)}}$.
$S^{(n)}$  has rational singularities, as the quotient of a smooth variety by a finite group (see \cite{Bu}, Proposition 4.1). 
Therefore its resolution $f_n:S^{[n]}\to S^{(n)}$ satifies   $R^i(f_n)_*\ko_{S^{[n]}}=0$.
\end{proof}


\section{The tangent bundle and tautological bundles}


For a smooth projective variety $X$ let $K(X)$ denote the Grothendieck group
generated by locally free sheaves, or equivalently, arbitrary coherent sheaves.
We will denote a sheaf and its class by the same symbol. $K(X)$ is endowed
with a ring structure which for classes of locally free sheaves $F$ and $G$
is given by $F\cdot G:=F\tensor G$, and a ring involution ${}\dual$, which 
extends $F\mapsto \khom(F,\ko_X)$ for a locally free sheaf. Thus for arbitrary 
coherent sheaves $F$, $G$ one has $F\dual=\sum_i(-1)^i\kext^i(F,\ko_X)$ and 
$F\dual\cdot  G=\sum_i(-1)^i\kext^i(F,G)$. For $f:X\to Y$ a morphism 
of smooth projective varieties there is a push-forward $f_!:K(X)\to K(Y);
\ F\mapsto \sum_i (-1)^iR^if_*(F)$, and a pullback $f^!$ which for 
a locally free sheaf $G$ on $Y$ is given by $f^!G=f^*G$.
In the paragraphs below, we will use the following well-known fact, 
which holds in greater
generality: let $G$ be an $S$-flat family of sheaves on $Y\times S$, and let
$p:Y\times S\to Y$ denote the projection. If $f:Y'\to Y$ is a morphism of
schemes, then $f^!(p_!G)=p_!(f_S^!G)$. Moreover, if 
the support of $G$ is finite over
$S$, there is an isomorphism of sheaves $f^*p_*G\isom p_*f_S^*G$.

In this section we describe recursion relations for the classes of the tangent
sheaf $T_n$ of $S^{[n]}$ and certain tautological sheaves $F^{[n]}$ on $S^{[n]}$
which are defined as follows: For any locally free sheaf $F$ on $S$ let
$F^{[n]}:=p_*(\ko_n\tensor q^*F)$ (with $p$ and $q$ as in diagram 
(\ref{dispA})). This extends to a group homomorphism
${}^{[n]}:K(S)\to K(S^{[n]})$. We will also write $p$ and $q$ for the 
projections of $S^{[n,n+1]}\times S$ to $S^{[n,n+1]}$ and $S$.

\begin{lemma}--- The following relation holds in $K(S^{[n,n+1]})$:
\beeq{tautorecursion}
\psi^!F^{[n+1]}=\phi^!F^{[n]}+\kl\cdot\rho^!F.
\eneq
\end{lemma}

\begin{proof}{} Let $F$ be a locally free sheaf on $S$. Apply the 
functor $p_*(\,.\,\tensor q^*F)$ to the exact sequence (\ref{basicseq2}) and
observe that
$$p_*(\phi_S^*(\ko_n)\tensor q^*F)=p_*(\phi_S^*(\ko_n\tensor q^*F))
=\phi^*(p_*(\ko_n\tensor q^*F))=\phi^*F^{[n]},$$
i.e. $p_!(\phi_S^!(\ko_n)\cdot q^!F)=\phi^!F^{[n]}$ and similarly
$p_!(\psi_S^!(\ko_{n+1})\cdot q^!F)=\psi^!F^{[n+1]}$.
Using $p\circ j=\id_{S^{[n,n+1]}}$ we get 
$$p_*(j_*(\kl)\tensor q^*F))=(p\circ j)_*(\kl\tensor (q\circ j)^*F)=
\kl\tensor \rho^*F,$$
i.e. $p_!(j_!(\kl)\cdot q^!F)) =\kl\cdot \rho^!F$.
\end{proof}

Next, we turn to the tangent sheaf:

\begin{proposition}\label{Tprop}--- The class of the tangent sheaf in 
$K(S^{[n]})$ is given by the relation
\beeq{}
T_n=\chi(\ko_S)\cdot 1- p_!(I_n\dual\cdot I_n).
\eneq
\end{proposition}

\begin{proof}{} We have the following isomorphisms
$$T_n\isom\Hom_p(I_n,\ko_n)\isom\Ext^1_p(\ko_n,\ko_n)\isom
p_*\kext^1_{\ko_{S^{[n]}\times S}}(\ko_n,\ko_n).$$
Here $\Ext_p$ are the higher derived functors of the composite functor
$p_*\circ\khom$.
For the first isomorphism see e.g.\ \cite{Lehn}. The last equality
is a consequence from the spectral sequence $R^ip_*\kext^j
\Rightarrow\Ext_p^{i+j}$
and the observation that the sheaves $\kext^*(\ko_n,\ko_n)$ are supported 
on the universal family $\Sigma_n$, which is finite over $S^{[n]}$ so that all
higher direct images vanish.
Moreover,
$$\kext^0(\ko_n,\ko_n)\isom\kext^0(\ko_{S^{[n]}\times S},\ko_n)\isom\ko_n$$
and, by $\kext^0(\ko_n,\ko_{S^{[n]}\times S})=0$, 
$\kext^1(\ko_n,\ko_{S^{[n]}\times S})=0$ also
$$\kext^2(\ko_n,\ko_n)\isom\kext^2(\ko_n,\ko_{S^{[n]}\times S})=\ko_n\dual.$$
Hence we get
\begin{eqnarray*}
T_n&=&p_*\kext^1(\ko_n,\ko_n)=p_!\ko_n-p_!(\ko_n\dual\cdot\ko_n)+p_!\ko_n\dual\\
&=&p_!(1)-p_!((1-\ko_n)\dual(1-\ko_n))=\chi(\ko_S)\cdot1-p_!(I_n\dual\cdot I_n)
\end{eqnarray*}
as an identity in $K(S^{[n]})$.
\end{proof}

\begin{proposition}--- The following relation holds in $K(S^{[n,n+1]})$:
\beeq{tangentrecursion}
\psi^!T_{n+1}=\phi^!T_n+\kl\cdot \sigma^! I_n\dual
+\kl\dual\cdot\sigma^!I_n\cdot\rho^!\omega_S\dual
-\rho^!(1-T_S+\omega_S\dual).
\eneq
\end{proposition}

\begin{proof}{} We have
\begin{eqnarray*}
\psi^!T_{n+1}&=&\psi^!(\chi(\ko_S)\cdot 1)
-\psi^!p_!(I_{n+1}\dual\cdot I_{n+1})\\
&=&\phi^!(\chi(\ko_S)\cdot 1)-p_!(\psi_S^!I_{n+1}\dual\cdot\psi_S^!I_{n+1})
\end{eqnarray*}
Now use (\ref{basicseq}) to replace $\psi_S^!I_{n+1}$ by 
$\phi_S^!I_n-p^!\kl\cdot\rho_S^!\ko_\Delta$, where $\Delta\subset S\times S$
denotes the diagonal. This yields
\begin{eqnarray*}
\psi^!T_{n+1}&=&\phi^!T_n-\rho^!p_!(\ko_\Delta\dual\cdot\ko_\Delta)\\
&&+p_!(\rho_S^!\ko_\Delta\cdot\phi_S^!I_n\dual)\cdot\kl
+p_!(\rho_S^!\ko_\Delta\dual\cdot\phi_S^!I_n).
\end{eqnarray*}
The two last summands can be simplified as follows:
$$p_!(\rho_S^!\ko_{\Delta}\cdot\phi_S^!I_n\dual)=
p_!(j_!(j^!\phi_S^!I_n\dual))=p_!j_!\sigma^!I_n\dual=\sigma^!I_n\dual$$
and similarly $p_!(\rho_S^!\ko_\Delta\dual\cdot\phi_S^!I_n)=
\sigma^!I_n\cdot\rho^!\omega_S\dual$,
since
$\rho_S^!\ko_\Delta\dual=\rho_S^!\Delta_!\omega_S\dual=j_!\rho^!\omega\dual$.
Finally, $p_!(\ko_\Delta\dual\cdot\ko_\Delta)=\ko_S-T_S+\omega_S\dual$,
which follows e.g.\ by Proposition \ref{Tprop}. \end{proof}

\section{The induction step}

We want to relate integrals on $S^{[n+1]}\times S^m$ to integrals on 
$S^{[n]}\times 
S^{m+1}$.

Let $Z=S^{[n,n+1]}\times S^{m}$. The maps $\phi$ and $\psi$ from diagram 
(\ref{dispB}) generalize to morphisms $\Psi=\psi\times\id_{S^m}:Z\to S^{[n+1]} 
\times S^m$ and $\Phi=\sigma\times\id_{S^m}:Z\to S^{[n]}\times S^{m+1}$.
For any $I\subset \{0,1,\ldots,m\}$ let $pr_I$ denote the projection from
$S^{[n+1]}\times S\times\ldots\times S$ to the product of the factors indexed
by $I$.

\begin{proposition}\label{technoprop}---
Let $f$ be a polynomial in the Chern classes of the following sheaves on 
$S^{[n+1]}\times S^m$:
$$pr_0^*T_{n+1},\,\,\, pr_{0i}^*I_{n+1},\,\,\, pr_{ij}^*\ko_\Delta,\,\,\, 
pr_i^*T_S\,\,\mbox{ for any }1\leq i,j\leq m.$$
Then there is a polynomial $\tilde f$, depending 
only on $f$, in the Chern classes of the analogously defined sheaves on 
$S^{[n]}\times S^{m+1}$ such that
$$\int_{S^{[n+1]}\times S^{m}}f=\int_{S^{[n]}\times S^{m+1}}\tilde f.$$
\end{proposition}

\begin{proof}{} The morphism $\Psi$ is generically finite of degree $n+1$,
 so that
$$\int_{S^{[n+1]}\times S^{m}}f=\frac{1}{n+1}\int_Z\Psi^*f.$$
Because of an index shift resulting from the insertion of the additional
factor $S$ between $S^{[n]}$ and $S^m$ we have
$$\Psi^!pr_i^*T_S=\Phi^!pr_{i+1}^*T_S,\,\Psi^!pr_{ij}^*\ko_\Delta=\Phi^!
pr_{i+1,j+1}^*\ko_\Delta.$$
Using (\ref{basicseq}) we get
$$\Psi^!pr_{0i}^*I_{n+1}=\Phi^!pr_{0,i+1}^*I_n+pr_Z^!\kl\cdot 
pr_{1,i+1}^*\ko_\Delta.$$
And finally by (\ref{tangentrecursion})
\begin{eqnarray}
\Psi^!pr_{0}^*T_{n+1}&=&\Phi^!pr_0^*T_n
+pr_{S^{[n,n+1]}}^!\kl\cdot pr_{01}^*I_n\dual\\
&&+pr_{S^{[n,n+1]}}^!\kl\dual\cdot pr_{01}^!I_n\cdot 
pr_1^!\omega_S\dual-pr_1^!(\ko_S-T_S+\omega_S\dual).
\end{eqnarray}
It follows that there are polynomials $f_\nu$ depending only on $f$, in the
Chern classes of the sheaves
$$pr_0^*T_{n},\, pr_{0i}^*I_{n},\, pr_{ij}^*\ko_\Delta,\, pr_i^*T_S$$
such that 
$$\int_{S^{[n+1]}\times S^{[m]}}f=\frac{1}{n+1}\int_Z\Psi^*f=
\int_Z\big(\sum_{\nu\geq0}\Phi^*f_\nu\cdot pr_Z^*(-c_1(\kl))^\nu\big).$$
As we are only trying to prove a general structure result 
we make no attempt to derive from the above 
recursion relations for the classes in the $K$-groups more explicit formulae
for the dependence of $f_\nu$ on $f$. 

Now, according to (\ref{pushdownofell}), the last integral equals:
$$\int_{S^{[n]}\times S^{m+1}}\Phi_*\big(\sum_{\nu\geq0}
\Phi^*f_\nu\cdot pr_Z^*(-c_1(\kl))^\nu\big)
=\int_{S^{[n]}\times S^{m+1}}\sum_{\nu\geq0}f_\nu\cdot c_\nu(-pr_{01}^*I_n).$$
The integrand in this expression is the polynomial $\tilde f$.
\end{proof}

\begin{proof}{ of Proposition \ref{propB}} Suppose we are given a 
polynomial $P$ in the Chern classes
 of $T_n$. Applying the proposition repeatedly, we may write
$$\int_{S^{[n]}}P=\int_{S^n} \tilde P$$
for some polynomial $\tilde P$, which depends only on $P$, in the Chern classes
of sheaves on $S^n$ of the form $pr_i^*T_S$ and $pr_{ij}^*\ko_\Delta$.
Any such expression $\int_{S^n}\tilde P$ can be universally reduced
to a polynomial expression of integrals of polynomials in the Chern classes
of $T_S$ (to see this for the Chern classes of $pr_{ij}^*\ko_\Delta$ 
one applies Riemann-Roch without denominators \cite{Jouanolou}).
 This finishes the proof.
\end{proof}


\section{A generalization}


We can generalize the methods used above to prove Proposition \ref{propB}
to cover integrals of polynomial expressions in the Chern classes  of 
tautological sheaves.
Let $F_1,\ldots,F_\ell\in K(S)$. We require that the rank $r_i$ of $F_i$
($i=1,\ldots,\ell$) is the same on all connected components of $S$.

\begin{theorem}\label{GeneralTheorem}--- Let $P$ be a polynomial
in the Chern classes of the tangent bundle $T_n$ of $S^{[n]}$ and the Chern
classes of $F_1^{[n]}$,\ldots,$F_\ell^{[n]}$.
Then there is a universal polynomial $\Ptilde$, depending only on P, in the 
Chern classes of $T_S$, the $r_1,\ldots,r_k$ and the Chern classes of   
$F_1$,\ldots,$F_\ell$ such that
$$\int_{S^{[n]}}P =\int_S
\Ptilde .$$
\end{theorem}

\begin{proof}{} The proof goes along similar lines as that of 
Proposition \ref{propB}. The result immediately follows from a modified version of 
Proposition \ref{technoprop}:
We now allow $f$ to be a polynomial in the $r_k$ and the Chern classes of 
the $$pr_0^*F_k^{[n+1]},pr_0^*T_{n+1},pr_{0i}^*I_{n+1}, pr_{ij}^*\ko_\Delta,
pr_i^*F_k,pr_i^*T_S$$ on $S^{[n+1]}\times S^{m}$, and get $\widetilde f$ 
to be a polynomial in the $r_k$ and the Chern classes of the 
analogously defined bundles on $S^{[n]}\times S^{m+1}$. To prove this use the 
recursion relation (\ref{tautorecursion}) and the formula 
$$c_i(\kl\otimes \rho^*F_k)=\sum_{a+b=i} \binom{r_k-b}{a}
c_1(\kl)^a\rho^*(c_b(F_k)),$$
to obtain
$$\int_{S^{[n+1]}\times S^{m}}f=
\int_Z\big(\sum_{\nu\geq0}\Phi^*f_\nu\cdot pr_Z^*(-c_1(\kl))^\nu\big),$$
where the $f_\nu$ are universal polynomials  in the $r_k$ and 
the Chern classes of 
$$pr_0^*F_k^{[n]},pr_0^*T_{n},pr_{0i}^*I_{n}, pr_{ij}^*\ko_\Delta,
pr_i^*F_k,pr_i^*T_S.$$ Then we again use  (\ref{pushdownofell}).
\end{proof}

We can use Theorem \ref{GeneralTheorem} to make predictions on the
algebraic structure of certain formulae: let $\Psi:K\lra H^\times$ be a 
multiplicative function, i.e.\ for any complex manifold $X$ we are given a
group homomorphism from the additive group $K(X)$ into the multiplicative
group of units in $H(X;\IQ)$, which is functorial with respect to pull-backs and
is polynomial in the Chern classes of its argument. 
The total Chern class,
the total Segre class or the Chern character of the determinant are of 
this type. Also let $\phi(x)\in\IQ[[x]]$ be a formal power series and, 
for a complex manifold $X$ of dimension $n$ put 
$\Phi(X):=\phi(x_1)\cdot\ldots\cdot\phi(x_n)\in H^*(X,\IQ)$ with
$x_1,\ldots,x_n$ the Chern roots of $T_X$.

Let $S$ be a smooth projective surface and let $x\in K(S)$. For any functions
$\Phi$ and $\Psi$ as above we define a power series in $\IQ[[z]]$ as follows: 
$$H_{\Psi,\Phi}(S,x):=\sum_{n=0}^\infty\int_{S^{[n]}}
\Psi(x^{[n]})\Phi(S^{[n]})\,z^n.$$

\begin{theorem}\label{uni}--- For each  integer $r$ 
there are universal power series $A_i\in\IQ[[z]]$, $i=1,\ldots,5$,
depending only on $\Psi$, $\Phi$ and $r$, such that for each $x\in K(S)$ of
rank $r$ (on every component of $S$) one has
$$H_{\Psi,\Phi}(S,x)=\exp(c_1^2(x)A_1+c_2(x)A_2+c_1(x)c_1(S)A_3+c_1^2(S)A_4+
c_2(S)A_5).$$
\end{theorem}

For simplicity we have suppressed the integrals $\int_S$ in the statement of 
the theorem and interpret the expressions $c_1(x)c_1(S)$ etc.\ as intersection 
numbers.

\begin{proof}{} Let $\kk_r:=\{(S,x)|S\mbox{ an algebraic surface}, x\in K(S),
{\rm rank}(x)=r\}$, and let $\gamma:\kk_r\to \IQ^5$ be the map
$(S,x)\mapsto(c_1^2(x),c_2(x),c_1(x)c_1(S),c_1^2(S),c_2(S))$. The images of 
the five elements 
$(\CP_2, r\cdot 1)$, $(\CP_2,\ko(1)+(r-1)\cdot 1$,
$(\CP_2,\ko(2)+(r-1)\cdot 1)$,
$(\CP_2,2\ko(1)+(r-2)\cdot 1)$, and
$(\CP_1^2,r\cdot 1)$ under $\gamma$ are the linearly independent vectors
$(0,0,0,9,3)$, $(1,0,3,9,3)$, $(4,0,6,9,3)$, $(4,1,6,9,3)$, and $(0,0,0,8,4)$. 

Now, if $S=S_1\sqcup S_2$ we may decompose $(S,x)\in \kk_r$
as $(S_1,x_1)+(S_2,x_2)$, where $x_i=x|_{S_i}$, and get
$\gamma(S,x)=\gamma(S_1,x_1)+\gamma(S_2,x_2)$. Moreover, there is
a decomposition of
the class of the tautological sheaf analogous to the decomposition 
(\ref{disjointunion}) of the Hilbert scheme:
$$x^{[n]}|_{S_1^{[n_1]}\times S_2^{[n_2]}}
=pr_1^*(x_1^{[n_1]})\cdot pr_2^*(x_2^{[n_2]}).$$
>From the multiplicative behaviour of $\Psi$ and $\Phi$ we deduce
$$\int_{S^{[n]}}\Psi(x^{[n]})\Phi(S^{[n]})=\sum_{n_1+n_2=n}
\int_{S_1^{[n_1]}}\Psi(x_1^{[n_1]})\Phi(S_1^{[n_1]})
\int_{S_2^{[n_2]}}\Psi(x_2^{[n_2]})\Phi(S_2^{[n_2]})$$
and get 
\beeq{multiplicprop}H_{\Psi,\Phi}(S,x)
=H_{\Psi,\Phi}(S_1,x_1)H_{\Psi,\Phi}(S_2,x_2).
\eneq
By Theorem \ref{GeneralTheorem}, the function $H_{\Psi,\Phi}:\kk_r\to\IQ[[z]]$
factors through $\gamma$ and a map $h:\IQ^5\lra\IQ[[z]]$.
As the image of $\gamma$ is Zariski dense in $\IQ^5$ we conclude from
(\ref{multiplicprop}), that
$$\log h(y_1+y_2)=\log h(y_1)+\log h(y_2),\mbox{ for all }y_1,y_2\in\IQ^5,$$
i.e.\ $\log h$ is a linear function. This proves the theorem.
\end{proof}


\section{Riemann-Roch numbers}


In this section we want to compute the Riemann-Roch numbers
$\chi(M)$ for line bundles and vector bundles on $S^{[n]}$.
Let  $L_n:=f^*g_*(\otimes_{i=1}^n pr_i^*L)^{\gothS_{n}}$, for $L$ a 
line bundle on 
$S$, where $f:S^{[n]}\to S^{(n)}$ and $g:S^n\to S^{(n)}$ are the  natural 
morphisms, and $pr_i:S^n\to S$ is the $i$-th projection. This gives a 
monomorphism $-_n:Pic(S)\to \Pic(S^{[n]}), L\mapsto L_n$.
It is well known that for $n\ge 2$
$$Pic(S^{[n]})=(Pic(S))_n\oplus \IZ E, \qquad E:=\det(\ko_S^{[n]}),$$
($E=\ko_S$ in case $n\le 1$).
Moreover, $c_1(E)=-\frac{1}{2}D$, where $D$ is the exceptional divisor
 of $S^{[n]}\to S^{(n)}$. Note that for $F\in K(S)$ 
$$\det(F^{[n]})=\det(F)_n\tensor E^{\rk(F)}.$$

\begin{lemma}\label{chivonLn}---
The Euler characteristics of $L_n\otimes E$, $L_n$ and $L$ are related 
by the formulae
$$\chi(L_{n})=\binom{\chi(L)+n-1}{n}\quad\mbox{ and }
\quad\chi(L_n\otimes E)=\binom{\chi(L)}{n}.$$
\end{lemma}

\begin{proof}{} Consider the cartesian diagram
\beeq{dispC}
\begin{array}{ccccc}
\widehat S^n_*&\stackrel{\widehat f}{\lra}&S^n_*&\\[1ex]
\scriptstyle{\widehat g}\Big\downarrow \phantom{\scriptstyle{\phi}}&&\phantom{\scriptstyle{\phi}}\Big\downarrow
\scriptstyle{g}\\[2ex]
S^{[n]}_*&\stackrel{ f}{\lra}&S^{(n)}_*,
\end{array}
\eneq
where $S^{(n)}_*$ is the open subscheme of zero cycles of length $n$ whose 
support consists of at least $(n-1)$ points and $S^{[n]}_*$ and 
$S^n_*$ are the preimages of $S^{(n)}_*$ 
under $f$ and $g$. 
It is easy to see that $\widehat S^n_*$ is the blow-up of $S^n_*$ along the
(disjoint) diagonals 
$$\Delta_{i,j}=\big\{(x_1,\ldots,x_n)\in S^n_*\bigm| x_i= x_j\big\}.$$
Let $\widehat\Delta_{i,j}$ denote the corresponding exceptional divisors, and
$\widehat\Delta:=\sum_{i<j} \widehat\Delta_{i,j}$ their sum.
The family $\Gamma\subset \widehat S^n_*\times S$ corresponding to the
classifying morphism $\widehat g$ is the union $\bigcup_{i=1}^n\Gamma_i$
of the graphs $\Gamma_i$ of the $i$-th projection $\widehat S^n_*\to 
S^n_*\to S$.
Projecting the short exact sequence
$$0\lra\ko_{\Gamma}\lra\bigoplus_{i=1}^n\ko_{\Gamma_i}
\lra\bigoplus_{i<j}\ko_{\Gamma_i\cap\Gamma_j}\lra0$$
to the base $\widehat S^n_*$ of these families we obtain a short exact
sequence
$$0\lra p_*\ko_{\Gamma}\lra\ko_{\widehat S^n_*}^{\oplus n}\lra \bigoplus_{i<j}\ko_{\hat\Delta_{ij}}=\ko_{\hat\Delta}\lra 0$$
of $\gothS_n$-linearized sheaves. $\gothS_n$ acts on the sheaf in the middle
by permutation of the factors.
If we take the determinant of the first 
homomorphism in this sequence, we obtain a short exact sequence
$$0\lra\det(p_*\ko_{\Gamma})\lra\ko_{\widehat S^n_*}^{\varepsilon}\lra \ko_{\widehat \Delta}\lra 0$$
of $\gothS_n$-linearized sheaves with equivariant homomorphisms, where the upper index $\varepsilon$ indicates
that the $\gothS_n$-linearization of $\ko_{\widehat S^n_*}$ is the standard
one twisted by the alternating character $\varepsilon:\gothS_n\to\IZ/2$.
(Recall that any two $G$-linearizations of a line bundle on a $G$-scheme differ 
by a character of the group $G$).
Thus we can identify
$\widehat g^*E=\det(p_*\ko_{\Gamma})$ with the $\gothS_n$-linearized subsheaf
$\ko_{\widehat S^n_*}(-\widehat\Delta)^{\varepsilon}\subset\ko_{\widehat S^n_*}^{\varepsilon}$, endowed with the alternating linearization.

Now $\widehat g^*L_n=\widehat f^*L^{\boxtimes n}$ and
$\widehat g^*(L_n\otimes E)=\widehat f^*L^{\boxtimes n}
\otimes \ko(-\widehat\Delta))$.
As $S^{[n]}$ is smooth and $S^{[n]}\setminus S^{[n]}_*$ has codimension $2$,
this gives 
\begin{eqnarray*}
H^0(S^{[n]},L_n)&=&H^0(S^{[n]}_*,L_n)=H^0(\widehat S^n_*,\widehat 
f^*L^{\boxtimes n})^{\gothS_n}\\
&=&H^0(S^n_*,L^{\boxtimes n})^{\gothS_n}=H^0(S^n,L^{\boxtimes n})^{\gothS_n}\\
&=&(H^0(S,L)^{\tensor n})^{\gothS_n}\isom S^nH^0(S,L),
\end{eqnarray*}
and, similarly,
\begin{eqnarray*}
H^0(S^{[n]},L_n\tensor E)&\isom& H^0(S^{[n]}_*,L_n\tensor E))\\
&=&H^0(\widehat S^{n}_*, \widehat f^*L^{\boxtimes n}
\otimes\ko(-\widehat\Delta))^{\gothS_n}\\
&=&H^0(S^n_*,L^\boxtimes)^{\gothS_n,\varepsilon}
=(H^0(S,L)^{\tensor n})^{\gothS_n,\varepsilon}\\
&\isom&\Lambda^n(H^0(S,L)).
\end{eqnarray*}
Here the notation $V^{\gothS_n,\varepsilon}$ means
$\{v\in V|\pi(v)=\varepsilon(\pi)\cdot v,\forall\,\pi\in\gothS_n\}$.
Taking dimensions, we get
\beeq{formulaforh0}
h^0(L_n)=\binom{h^0(L)+n-1}{n}\quad\mbox{and}
\quad h^0(L_n\tensor E)=\binom{h^0(L)}{n}.
\eneq

Now let $H$ be an ample line bundle on $S$. By Grothendieck's construction of
the Hilbert scheme it follows that $H_n\tensor E\isom\det(H^{[n]})$
is very ample for sufficiently ample $H$ on $S$. Applying this to 
$L\tensor H^k$ for
sufficiently large $k$, We conclude that
$$\chi(L_n\tensor H^k_n\tensor E)
=h^0(L_n\tensor H_n^k\tensor E)=\binom{h^0(L\tensor H^k)}{n}
=\binom{\chi(L\tensor H^k)}{n}.
$$
Evaluating this equation of polynomials in $k$ at $k=0$ one
finds $\chi(L\tensor E)=\binom{\chi(L)}{n}$ for all line bundles $L$.

Let $L$ be an arbitrary line bundle on $S$ and let $H$ be ample $S$. 
Then $H^{\boxtimes n}$ is an $\gothS_n$-linearized ample line bundle on $S^{n}$ and 
descends to an ample line bundle $H^{(n)}$ on $S^{(n)}$. For sufficiently 
large $k$, the line bundles $L\tensor H^k$ on $S$ and $(L\tensor H^k\tensor
 \omega_S\dual)^{(n)}$ on $S^{[n]}$ will be very ample. In particular, $f^*(L\tensor H^k\tensor\omega_S\dual)^{(n)}\isom
L_n\tensor H_n^k\tensor\omega_{S^{[n]}}\dual$ is globally generated and
big. It follows from the Grauert-Riemenschneider vanishing theorem
\cite{GrauertRiemenschneider} that
$H^i(S^{[n]},L_n\tensor H_n^k)=0$ for all $i>0$. This gives
$$\chi(L_n\tensor H_n^k)=h^0(L_n\tensor H_n^k)
=\binom{h^0(L\tensor H^k)+n-1}{n}=
\binom{\chi(L\tensor H^k)+n-1}{n}$$
for all sufficiently large $k$. As both sides are polynomials in $k$, we may 
take $k=0$.
\end{proof}

Consider the following power series in $\IQ[a,y][[z]]$:
$$f_{y,a}:=\sum_{n\ge 0}\binom{y-a(n-1)}{n}z^n,
\quad g_{y,a}:=\sum_{n\ge 0}\frac{y}{y-an}\binom{y-an}{n}z^n.$$
$f_{y,a}$ and $g_{y,a}$ play a role in combinatorics:
Take $k$ points on a line with a distance of $1$ from one to the next.
Then $\binom{k-a(n-1)}{n}$ is the number of ways to choose $n$ of them
with a distance $>a$ among each other, and 
$\frac{k}{k-an}\binom{k-an}{n}$ is the same number if the $k$ points lie on
a circle \cite{Rio}. 

In the proof of the following lemma we make use of the 
B\"urmann-Lagrange expansion
formula in the following form: if $L(z)$ is a power series, and $\alpha(z)$ is
a power series with $\alpha(0)=\alpha'(0)=0$, then
$$\sum_{n\geq 0}\frac{1}{n!}\frac{d^n}{dz^n}(L(z)\alpha(z)^n)
=\frac{L(\zeta)}{1-\alpha'(\zeta)}\Big|_{\zeta=z-\alpha(\zeta)}.$$
For a proof in an analytic context see \cite[p 128-135]{WhitWat}.

\begin{lemma}\label{powseries}--- The power series $f_{y,a}$ and $g_{y,a}$ 
are related as follows:
$$g'_{y,a}=yf_{y-2a-1,a},\qquad g_{y,a}=g_{1,a}^y,\qquad 
f_{y,a}=g_{1,a}^y\cdot f_{0,a}.$$
\end{lemma}

\begin{proof}{} Using new formal variables $u,x$ and $v$ that are related by
$z=u^a$, $u=x+x^{a+1}$, and $v=x^a$, so that $z=v(1+v)^a$, and, in the last 
line of the computation, the B\"urmann-Lagrange formula, we have
\begin{eqnarray*}
f_{y,a}(z)&=&\sum_{n\geq 0}\binom{y-a(n-1)}{n}z^n\\
&=&\sum_{n\geq 0}\binom{(a+1)n-(y+a+1)}{n}(-z)^n\\
&=&\sum_{n\geq 0}u^{y+a+1}\frac{1}{n!}\frac{d^n}{du^n}
\frac{(-u^{a+1})^n}{u^{y+a+1}}\\
&=&\left(\frac{u}{x}\right)^{y+a+1}\frac{1}{1+(a+1)x^a}=
(1+v)^y\cdot\frac{(1+v)^{a+1}}{1+(a+1)v}
\end{eqnarray*}
This gives $f_{0,a}=(1+v)^{a+1}/(1+(a+1)v)$ and $f_{y,a}(z)
=f_{0,a}(z)\cdot (1+v(z))^{y}$.
Now $dz=(1+(a+1)v)(1+v)^{a-1}dv$ and hence
$$\frac{d}{dz}
(1+v(z))=\frac{1}{(1+v)^{a-1}(1+(a+1)v)}=\frac{f_{0,a}}{(1+v)^{2a}}.$$
Differentiating $g_{y,a}$ we find
\begin{eqnarray*}
\frac{d}{dz}g_{y,a}(z)&=&\sum_{n\geq 1}\frac{yn}{y-an}\binom{y-an}{n}z^{n-1}\\
&=&\sum_{n\geq 1}y\binom{y-an-1}{n-1}z^{n-1}\\
&=&yf_{y-2a-1,a}=y(1+v)^{y-1-2a}f_{0,a}\\
&=&y(1+v)^{y-1}\cdot\frac{f_{0,a}}{(1+v)^{2a}}=\frac{d}{dz}(1+v)^y.
\end{eqnarray*}
As both $g_{y,a}$ and $(1+v)^y$ are power series in $z$ with constant term 1, we
have $g_{y,a}=(1+v)^y$. Collecting the proven relations we conclude that
$g_{1,a}=1+v$, $g_{y,a}=g_{1,a}^y$ and $f_{y,a}=g_{1,a}^yf_{0,a}$.
\end{proof}

\begin{theorem}\label{chithm}--- Let $K$ denote  the canonical divisor of $S$.
For each $r\in \IZ$, there exist universal power series $A_r,B_r\in \IQ[[z]]$, 
such that for all $L\in Pic(S)$
$$\sum_{n\ge 0} \chi(L_n\tensor E^r))z^n=g_{1,r^2-1}(z)^{\chi(L)}\cdot
f_{0,r^2-1}(z)^{\frac{\chi(\ko_S)}{2}}\cdot  A_r(z)^{KL-\frac{K^2}{2}}
\cdot B_r(z)^{K^2}.$$
Moreover, $A$ and $B$ satisfy the symmetry relations $A_{-r}=1/A_r$  
and $B_{-r}=B_r$ for arbitrary $r$, and $A_r=B_r=1$ for $r=-1,0,1$.
In particular, 
$$\chi(L_n\tensor E^{\pm 1})= \binom{\chi(L)}{n}.$$ 
If $S$ is a K3-surface, then 
$$\chi(L_n\tensor E^r)=\binom{\chi(L)-(r^2-1)(n-1)}{n},$$
and if 
$S$ is an abelian surface, then
$$\chi(L_n\tensor E^r)=\frac{\chi(L)}{n}\binom{\chi(L)-(r^2-1)n-1}{n-1}.$$
\end{theorem}

\begin{proof}{} Let $F=L+r\cdot 1\in K(S)$. Then $\det(F^{[n]})=
L_n\tensor E^r$,
$c_1(F)=c_1(L)$, $c_2(F)=0$. Therefore, by Theorem \ref{uni},   
\begin{eqnarray*}
\sum_{n\geq 0}\chi(L_n\tensor E^r)z^n
&=&Z_1(z)^{K^2}Z_2(z)^{KL}Z_3(z)^{c_2(S)}Z_4(z)^{L^2}\\
&=&A_r(z)^{K L-\frac{K^2}{2}}B_r(z)^{K^2}
F_r(z)^{\frac{\chi(\ko_S)}{2}}G_r(z)^{\chi(L)},
\end{eqnarray*}
 for suitable power series $Z_i,A_r,B_r,G_r,F_r\in \IQ[[z]]$.
For the second equality we have used the identities
$\chi(L)=\frac{L(L-K)}{2}+\chi(\ko_S)$ and $\chi(\ko_S)=
\frac{c_1(S)^2+c_2(S)}{12}$.

It is well-known that $\omega_{S^{[n]}}=(\omega_S)_n$. We get by 
Serre duality
$$\chi(L_n\tensor E^{-r}))=\chi(\omega_{S^{[n]}}\tensor L_n\dual\tensor E^r)
=\chi((\omega_S\tensor L\dual)_n\tensor E^r)).$$
Using $\chi(L)=\chi(\omega_S\tensor L\dual)$ and $K(K-L)-\frac{K^2}{2}=
-(KL-\frac{K^2}{2})$, this 
gives $B_{-r}=B_r$, $F_{-r}=F_r$, $G_{-r}=G_r$ and $A_{-r}=1/A_{r}$.

To  determine $F_r$ and $G_r$ explicitly, let $S$ be a K3-surface. 
Then by \cite{Bea} the Hilbert scheme $X:=S^{[n]}$ is an irreducible  
symplectic complex manifold. 

There exists a natural quadratic form  
$q$ on $H^2(X,\IZ)$ (see \cite{Bea},\cite{Fu},\cite{Huybas}), which on
$H^{1,1}(X)$ is defined as follows: let $\omega\in H^2(X,\IZ)$ be
the  everywhere non-degenerate holomorphic $2$-form,
normalized by $\int_X \omega^n\bar\omega^n =n!$. Then for $\alpha\in
H^{1,1}(X)$:
$$q(\alpha)=\frac{1}{(n-1)!}\int_X(\omega\bar\omega)^{n-1}\alpha^2.$$
Moreover, there exists a universal polynomial 
$h\in \IQ[z]$ such that $\chi(M)=h(q(c_1(M)))$ for all $M$ in $Pic(X)$.
For $L\in Pic(S)$ Beauville \cite[Lemma 9.1, 9.2]{Bea}  showed that 
$q(c_1(L_n))=L^2$, $q(c_1(E))=-2(n-1)$, 
and $E$ and $L_n$ are orthogonal with  respect to $q$; 
thus $q(c_1(L_n\otimes E^r))=L^2-2r^2(n-1)$.

The polynomial $h$ is determined by the formula of Lemma \ref{chivonLn} 
$$\chi(L_n)=\binom{\chi(L)+n-1}{n}=\binom{L^2/2+n+1}{n}=\binom{q(L)/2+n+1}{n}
$$
applied to sufficiently many $L$ on $S$ with distinct values of $L^2$:
\begin{eqnarray*}\chi(L_n\tensor E^r)&=&h(q(c_1(L_n\tensor E^r)))  \\[1ex]
&=&\binom{q(c_1(L_n\tensor E^r))/2+n+1}{n}\\[1ex]
&=&\binom{\chi(L)-(r^2-1)(n-1)}{n}.
\end{eqnarray*}
As $\chi(\ko_S)=2$, we get 
$F_r\cdot G_r^{\chi(L)}=f_{\chi(L),r^2-1}$. It follows from
Lemma \ref{powseries} that we can identify $F_r=f_{0,r^2-1}$ and
$G_r= g_{1,{r^2-1}}$. Finally, the computations of Lemma \ref{chivonLn} show
that the equation $\chi(L\tensor E^r)=\binom{\chi(L)-(r^2-1)(n-1)}{n}$
holds for all surfaces, if we restrict to small ranks $r=-1,0,1$. This 
implies
that $A_r=B_r=1$ for $r=-1,0,1$.
\end{proof}

\begin{remark}\label{maple}--- Using Bott's residue formula the
coefficients of $A_r$ and $B_r$ can be determined.
We computed $A_r$ and $B_r$ up to order $8$
in $z$. 
Up to order $5$ in $z$ we get
%
%
%
%
%
\begin{eqnarray*}
\log(A_r)&\equiv&\left(\frac{1}{6}r-\frac{1}{6}{r}^{3}\right)z^2
+\left(\frac{1}{5}r-\frac{5}{8}r^3+\frac {17}{40}r^5\right)z^3\\
&&
+\left(\frac{29}{140}r-\frac{209}{180}r^3+\frac{88}{45}r^5
-\frac{631}{630}r^7\right)z^4\\
&&
+\left(\frac{13}{63}r-\frac{31259}{18144}r^3+\frac{16979}{3456}r^5
      -\frac{69619}{12096}r^7+\frac{171215}{72576}r^9\right)z^5
\end{eqnarray*}
\begin{eqnarray*}
B_r&\equiv&1+\left (\frac{1}{24}{r}^{2}-\frac{1}{24}{r}^{4}\right ){z}^{2}
+\left ({\frac {29}{360}}{r}^{2}-{\frac
{31}{144}}{r}^{4}+{\frac {97}{720}}{r}^{6}
\right){z}^{3}\\
&&+\left(\frac{139}{1260}r^2-\frac{3053}{5760}r^4+
\frac{2273}{2880}r^6-\frac{14899}{40320}r^8\right)z^4\\
&&+\left(\frac{187}{1400}r^2-\frac{6257}{6480}r^4
+\frac{421267}{172800}r^6-\frac{311701}{120960}r^8
+\frac{503377}{518400}r^{10}\right ){z}^{5}
\end{eqnarray*}
The method for showing this result is as follows. By  Theorem \ref{chithm}
it is enough to compute $\chi((kH)_n+rE)$ for $S=\IC\IP_2$ and $H$ the hyperplane
bundle.
By the Riemann-Roch theorem we have to compute
$$\int_{\IC\IP_2^{[n]}}td_{\IC\IP_2^{[n]}} \cdot \exp((kH)_n+rE).$$
For this we use Bott's residue formula. The maximal torus $\Gamma$ of 
$SL(2,\IC)$ acts on   $\IC\IP_2^{[n]}$ with finitely many fix points $Z$, which
together with the structure of the tangent space as $\Gamma$-module are
explicitly described in \cite{GE-SAS}. In \cite{GE-SAS2} also the 
structure of the fibres $(nH)^{[n]}$ as $\Gamma$ module is given.
This information is enough to apply Bott's residue formula to compute
the above intersection number following the strategy outlined in 
\cite{GE-SAS2}. The computations are carried out with a 
suitable Maple program. 
\end{remark}

\begin{remark}\label{maple2}
Theorem \ref{MainTheorem} allows us to use the Bott residue formula 
to compute the Chern numbers of $S^{[n]}$ for all $n$  and  any surface $S$:
We computed them for $n\le 7$. As an illustration we give the numbers for
$S$ a $K3$ surface (in this case the odd Chern classes vanish)
and $n\le 4$. For a partition $(n_1,\ldots,n_r)$ of $2n$ we write
$(n_1,\ldots,n_r)$ for $c_{n_1}(S^{[n]})\cdot\ldots\cdot c_{n_r}(S^{[n]})$. We obtain
$(4) = 324$, 
$(2^2) = 828$,
$(6) = 3200$, 
$(4,2) = 14720$
$(2^3)= 36800$,
$(8) = 25650$,
$(6,2) = 182340$,
$(4^2) = 332730$,
$(4,2^2) = 813240$,
$(2^4)= 1992240$.
 
The strategy is similar to that of Remark \ref{maple}: by Theorem  \ref{MainTheorem}
it is enough to compute the Chern numbers in case $S=\IP_2$ and $S=\IP_1\times
\IP_1$. In both  cases we have an action of a 2-dimensional torus $\Gamma$
on $S^{[n]}$ with finitely many fix points. We proceed as above, again
using  a suitable Maple program.

It is remarkable that for any surface $S$ and all $n\le 7$ all the 
Chern numbers of $S^{[n]}$ are polynomials in $c_1(S)^2$ and $c_2(S)$ 
with nonnegative coefficients (this was observed by G. Thompson).
 G. H\"ohn has used these numbers to check the conjecture
of \cite{DMVV} about the elliptic genus of the Hilbert scheme of points $S^{[n]}$
on a $K3$ surface  for $n\le 6$.
\end{remark}

Finally, we compute the holomorphic Euler characteristics $\chi(F^{[n]})$
of the tautological bundles $F^{[n]}$ on $S^{[n]}$. Related results
are shown, using different methods in the recent preprint \cite{D}.

\begin{proposition}\label{chifn}--- i) 
Let $F$ be a vector bundle on $S$. Then 
$$
\sum_{n,i} h^i(S^{[n]},F^{[n]})z^it^{n-1}=\Big(\sum_j h^j(S,F)z^j\Big)
\frac{(1+zt)^{h^1(\ko_S)}}{(1-t)^{h^0(\ko_S)}(1-z^2t)^{h^2(\ko_S)}}.
$$
In particular, if $S$ is connected, then $h^0(S^{[n]},F^{[n]})=h^0(S,F)$,
and if furthermore $h^i(S,\ko_S)=0$ for $i>0$, then 
$h^i(S^{[n]},F^{[n]})=h^i(S,F)$ for all $i$. 

ii) For all $F\in K(S)$ one has
$$\chi(F^{[n]})=\chi(F)\binom{\chi(\ko_S)+n-2}{n-1}.$$
\end{proposition}

\begin{proof}{} Part ii) follows from i) by putting $z=-1$. 
In order to prove part i) consider the cartesian diagram for $n\geq 1${}
$$\begin{array}{ccccc}
\Sigma_n&\stackrel{\overline f}{\lra}&Z_n&\\[1ex]
\scriptstyle{p}\Big\downarrow \phantom{\scriptstyle{p}}&&\phantom{\scriptstyle{\overline p}}\Big\downarrow
\scriptstyle{\overline p}\\[2ex]
S^{[n]}&\stackrel{ f}{\lra}&S^{(n)},
\end{array}$$
where $f$ is the Hilbert-Chow morphism and $$Z_n:=\big\{(\eta,x)\in S^{(n)}\times S\bigm| x\in supp(\eta)\big\}$$
is the image of the universal family 
$\Sigma_n\subset S^{[n]}\times S$ in $S^{(n)}\times S$.
We denote by $q:\Sigma_n\to S$, $\overline q:Z_n\to S$ the restriction of the 
projection. It is easy to see that there is an isomorphism $Z_n\isom S^{(n-1)}\times S$
which identifies $\overline q$ with the projection to the second factor.
The maps $p$ and $\overline p$ are finite. By proposition \ref{hilbertchowvanishing}
the higher direct images $R^if_*\ko_{S^{[n]}}$ vanish. 
An easy spectral sequence argument shows that then the higher direct image
sheaves $R^i\overline f_*\ko_{\Sigma_n}$ must vanish as well. Moreover, 
$\overline f$ is a birational morphism of integral varities and $Z_n$ is
normal. This shows that $\overline f_*\ko_{\Sigma_n}=\ko_{Z_n}$. In particular,
for any locally free sheaf $G$ on $Z_n$ one has $H^i(\Sigma_{n},\overline f^*G)=H^i(Z_n,G)$.

By definition, $F^{[n]}=p_*q^*F$. As $p$ is finite, we get
\begin{eqnarray*}
H^i(S^{[n]},F^{[n]})&\isom& H^i(\Sigma_n,q^* F)=H^i(\Sigma_n,\overline f^*\overline q^*F)\\
&=&H^i(Z_n,\overline q^*F)=H^i(S^{(n-1)}\times S,\ko_{S^{(n-1)}}\boxtimes F)
\end{eqnarray*}
for all $i$. By the K\"unneth formula this gives
$$\sum_i h^i(S^{[n]},F^{[n]})z^i=\Big(\sum_ih^i(S,F)z^i\Big)\cdot
\Big(\sum_ih^i(S^{(n-1)},\ko_{S^{(n-1)}})z^i\Big).$$
By proposition \ref{hilbertchowvanishing}, $h^i(S^{(n-1)},\ko_{S^{[n-1]}})
=h^i(S^{[n-1]},\ko_{S^{[n-1]}})$. Hence
$$\sum_{n}\sum_{i}h^i(S^{[n]},F^{[n]})z^it^{n-1}
=\Big(\sum_ih^i(S,F)z^i\Big)\cdot\Big(\sum_{\nu}
\sum_ih^i(S^{[\nu]},\ko_{S^{[\nu]}})z^it^\nu\Big)
$$
The second factor on the right hand side was computed in \cite[Proposition 3.3]{Goe} and equals
$$\frac{(1+zt)^{h^1(S,\ko_S)}}{(1-t)^{h^0(S,\ko_S)}(1-z^2t)^{h^2(S,\ko_S)}}.$$
\end{proof}


\begin{small}
\noindent
Geir Ellingsrud\\
Matematisk Institutt, Universitetet i Oslo\\
Boks 1053 Blindern, N-0136 Oslo, Norway\\
e-mail: \verb{ellingsr@math.uio.no{

\medskip

\noindent
Lothar G\"ottsche\\
Abdus Salam International Center for Theoretical Physics\\
P.O.~Box 586, 34100 Trieste, Italy\\
e-mail: \verb{gottsche@ictp.trieste.it{

\medskip

\noindent
Manfred Lehn\\
Mathematisches Institut der Georg-August-Universit\"at\\
Bunsenstra\ss{}e 3-5, D-37073 G\"ottingen, Germany\\
e-mail: \verb{lehn@uni-math.gwdg.de{

\end{small}

\end{document}